\documentclass[a4paper,12pt]{article}
\usepackage{amsmath}
\usepackage{latexsym}
\numberwithin{equation}{section}

\def\R{{\bf R}}

\def\d{\displaystyle}
\def\e{{\varepsilon}}

\def\wt{\widetilde}
\def\wh{\widehat}

\def\v#1{\mbox{\boldmath $#1$}}

\newtheorem{thm}{Theorem}[section]

\newtheorem{lem}{Lemma}[section]
\newtheorem{prop}{Proposition}[section]
\newtheorem{rem}{Remark}[section]

\title{The sharp upper bound of the lifespan of solutions
to critical semilinear wave equations
in high dimensions
\footnote{
This paper already appeared as
Hokkaido University Preprint Series in Mathematics $\#$969 (24 September 2010)
http://eprints3.math.sci.hokudai.ac.jp/2118/
}
}
\author{
Hiroyuki Takamura\footnote{Department of Complex and Intelligent Systems,
Faculty of Systems Information Science,
Future University Hakodate,
116-2 Kamedanakano-cho,
Hakodate, Hokkaido 041-8655, Japan.
e-mail : takamura@fun.ac.jp}
\quad
and
\quad
Kyouhei Wakasa\footnote{The 4th year of the undergraduate,
Department of Complex Systems, School of Systems Information Science,
Future University Hakodate,
116-2 Kamedanakano-cho, Hakodate, Hokkaido 041-8655, Japan.
e-mail : b1007206@fun.ac.jp}
}
\date{
\[
\begin{array}{l}
\mbox{\scriptsize{\bf Keywords:} lifespan, semilinear wave equation, critical exponent, high dimensions}\\
\mbox{\scriptsize{\bf MSC2010:} Primary 35L71; Secondary 35B33, 35B44, 35E15}
\end{array}
\]
}
\begin{document}
\maketitle
\begin{abstract}
\par
The final open part of Strauss' conjecture on semilinear wave equations
was the blow-up theorem for the critical case in high dimensions.
This problem was solved by Yordanov and Zhang \cite{YZ06},
or Zhou \cite{Z07} independently.
But the estimate for the lifespan, the maximal existence time,
of solutions was not clarified in both papers.
\par
In this paper, we refine their theorems and introduce a new iteration argument
to get the sharp upper bound of the lifespan.
As a result, with the sharp lower bound by Li and Zhou \cite{LZ95},
the lifespan $T(\e)$ of solutions of
$u_{tt}-\Delta u=u^2$ in $\R^4\times[0,\infty)$
with the initial data $u(x,0)=\e f(x),u_t(x,0)=\e g(x)$
of a small parameter $\e>0$,
compactly supported smooth functions $f$ and $g$, has an estimate
\[
\exp\left(c\e^{-2}\right)\le T(\e)\le\exp\left(C\e^{-2}\right),
\]
where $c$ and $C$ are positive constants depending only on $f$ and $g$.
This upper bound has been known to be the last open optimality of the general theory for
fully nonlinear wave equations.
\end{abstract}

\section{Introduction}
\par
First we shall outline the general theory on the initial value problem
for fully nonlinear wave equations,
\begin{equation}
\label{GIVP}
\left\{
\begin{array}{l}
u_{tt}-\Delta u=H(u,Du,D_xDu) \quad \mbox{in}\quad\R^n\times[0,\infty),\\
u(x,0)=\e f(x),\ u_t(x,0)=\e g(x),
\end{array}
\right.
\end{equation}
where $u=u(x,t)$ is a scalar unknown function of space-time variables,
\[
\begin{array}{l}
Du=(u_{x_0},u_{x_1},\cdots,u_{x_n}),\ x_0=t,\\
D_xDu=(u_{x_ix_j},\ i,j=0,1,\cdots,n,\ i+j\ge1),
\end{array}
\]
$f,g\in C^\infty_0(\R^n)$ and $\e>0$ is a small parameter.
Let
\[
\wh{\lambda}=(\lambda;\ (\lambda_i),i=0,1,\cdots,n;
\ (\lambda_{ij}),i,j=0,1,\cdots,n,\ i+j\ge1). 
\]
Suppose that the nonlinear term $H=H(\wh{\lambda})$ is a sufficiently smooth function with
\[
H(\wh{\lambda})=O(|\wh{\lambda}|^{1+\alpha})
\]
in a neighbourhood of $\wh{\lambda}=0$, where $\alpha\ge1$ is an integer.
Let us define the lifespan $\wt{T}(\e)$ by
\[
\wt{T}(\e)=\sup\{t>0\ :\ \exists\mbox{solution $u(x,t)$ of (\ref{GIVP})
for arbitrarily fixed $(f,g)$.}\}.
\]
When $\wt{T}(\e)=\infty$, the problem (\ref{GIVP}) admits
a global in time classical solution,
while we only have a local in time solution on $t\in[0,\wt{T}(\e))$
when $\wt{T}(\e)<\infty$.
For local in time solutions, one can measure the global stability of a zero solution
by orders of $\e$.
Because the uniqueness of the solution of (\ref{GIVP})
may yield that $\d\lim_{\e\rightarrow0}\wt{T}(\e)=\infty$.
Such a uniqueness theorem can be found in Appendix of John \cite{J90} for example.
For $n=1$, we have no time decay of solutions even for the free case,
so that there is no possibility to obtain any global in time solution
of (\ref{GIVP}).
In this paper we assume $n\ge2$ for the simplicity. 

In Chapter 2 of Li and Chen \cite{LC92},
we have long histories on the estimate for $\wt{T}(\e)$.
The lower bounds of $\wt{T}(\e)$ are summarized in the following table.
Let $a=a(\e)$ satisfy
\begin{equation}
\label{a}
a^2\e^2\log(a+1)=1
\end{equation}
and $c$ stand for a positive constant independent of $\e$.
Then,
due to the fact that it is impossible to obtain an $L^2$ estimate for $u$ itself
by standard energy methods, we have
\begin{center}
\begin{tabular}{|c||c|c|c|}
\hline
$\wt{T}(\e)\ge$  & $\alpha=1$ & $\alpha=2$ & $\alpha\ge3$\\
\hline
\hline
$n=2$ &
$\begin{array}{l}
ca(\e)\\
\quad\mbox{in general case},\\
c\e^{-1}\\
\quad\mbox{if}\ \int_{\R^2}g(x)dx=0,\\
c\e^{-2}\\
\quad\mbox{if}\ \partial^2_uH(0)=0
\end{array}
$
&
$\begin{array}{l}
c\e^{-6}\\
\quad\mbox{in general case},\\
\exp(c\e^{-2})\\
\quad\mbox{if}\ \partial^b_uH(0)=0\ (b=3,4)
\end{array}
$
& $\infty$ \\
\hline
$n=3$ &
$\begin{array}{l}
c\e^{-2}\\
\quad\mbox{in general case},\\
\exp(c\e^{-1})\\
\quad\mbox{if}\ \partial^2_uH(0)=0
\end{array}
$
& $\infty$ & $\infty$ \\\hline
$n=4$ &
$\begin{array}{l}
\exp(c\e^{-2})\\
\quad\mbox{in general case},\\
\infty\\
\quad\mbox{if}\ \partial^2_uH(0)=0
\end{array}$
& $\infty$ & $\infty$ \\
\hline
$n\ge5$ & $\infty$ & $\infty$ & $\infty$\\
\hline
\end{tabular} 
\end{center}
We note that the lower bound in the case where $n=4$ and $\alpha=1$ is $\exp(c\e^{-1})$ 
in general case in Li and Chen \cite{LC92}.
But later, Li and Zhou \cite{LZ95} improve this part.
The remarkable fact is that {\bf all these lower bounds are known to be sharp
except for $\v{n=4}$ and $\v{\alpha=1}$}.
See Li and Chen \cite{LC92} for references on the whole history.
\par
Our purpose in this paper is to show this remained sharpness of the lower bound
by giving a sharp blow-up theorem for $u_{tt}-\Delta u=u^2$ in $\R^4\times[0,\infty)$. 
Including this situation,
we consider the initial value problem for semilinear wave equations of the form,
\begin{equation}
\label{IVP}
\left\{
\begin{array}{l}
u_{tt}-\Delta u=|u|^p \quad \mbox{in}\quad\R^n\times[0,\infty),\\
u(x,0)=\e f(x),\ u_t(x,0)=\e g(x),
\end{array}
\right.
\end{equation}
where $p>1$.
Let us define the lifespan $T(\e)$ by
\[
T(\e)=\sup\{t>0\ :\ \exists\mbox{solution $u(x,t)$ of (\ref{IVP})
for arbitrarily fixed $(f,g)$.}\},
\]
where \lq\lq solution" means the classical one if $p\ge2$, or the
weak one which is the solution of associated integral equations to (\ref{IVP}) if $1<p<2$.
Then we have the following Strauss' conjecture.
There exists a critical number $p_0(n)$ such that
\[
\begin{array}{lll}
T(\e)=\infty & \mbox{if $p>p_0(n)$ and $\e$ is \lq\lq small"}
& \mbox{(global in time existence)},\\
T(\e)<\infty & \mbox{if $1<p\le p_0(n)$}
& \mbox{(blow-up in finite time)}.
\end{array}
\]
As in Section 4 in Strauss \cite{St89},
$p_0(n)$ is a positive root of the quadratic equation
\begin{equation}
\label{gamma}
\gamma(p,n)\equiv 2+(n+1)p-(n-1)p^2=0.
\end{equation}
That is,
\begin{equation}
\label{p_0(n)}
p_0(n)=\frac{n+1+\sqrt{n^2+10n-7}}{2(n-1)}
\end{equation}
and one should remark that $p_0(4)=2$.
This number comes from the integrability of a weight function
$(1+|t-|x||)^{(n-1)p/2-(n+1)/2}$ in the iteration argument.
Such a weight function arises from the space-time integration
of $(1+t+|x|)^{(n-1)/2}$ which is a decay of a solution to free wave equation.
Note that we have another story for non-compactly supported data,
such as $T(\e)<\infty$ even for the supercritical case $p>p_0(n)$
if the spatial decay at infinity of the data is weak.
All the results in this direction
are summarized in Takamura, Uesaka and Wakasa \cite{TUW10}.
\par
Strauss' conjecture was first verified by John \cite{J79} for $n=3$
except for $p=p_0(3)$. Later, Glassey \cite{G81a, G81b} verified
this for $n=2$ except for $p=p_0(2)$.
Both critical cases were studied by Schaeffer \cite{Sc85}.
In high dimensions, $n\ge4$, the subcritical case was proved by Sideris \cite{Si84}.
For the supercritical case, there were many partial results.
The final result was given by Georgiev, Lindblad and Sogge \cite{GLS97}.
The critical case in high dimensions was obtained
by Yordanov and Zhang \cite{YZ06}, or Zhou \cite{Z07} independently.
In this way, the open part of the conjecture has been disappeared.
\par
For (\ref{IVP}), we have precise results on bounds of the lifespan 
in low dimensions, $n=2,3$, by virtue of the positivity of the fundamental solution.
Actually we know that
\begin{equation}
\label{lifespan1}
\lim_{\e\rightarrow0}\e^{2p(p-1)/\gamma(p,n)}T(\e)>0
\quad\mbox{exists for}\ l(n)<p<p_0(n),
\end{equation}
where $l(3)=1$ and $l(2)=2$.
This result was proved by Lindblad \cite{L90} for $n=3$ 
and by Zhou \cite{Z93} for $n=2$.
In Lindblad \cite{L90}, it was also proved that
for $n=2$ and $p=2$ we have
\begin{equation}
\label{lifespan2}
\begin{array}{ll}
\d\lim_{\e\rightarrow0}a(\e)^{-1}T(\e)>0 & \d\mbox{exists if}\ \int_{\R^2}\!\!g(x)dx\neq0,\\
\d\lim_{\e\rightarrow0}\e T(\e)>0 & \d\mbox{exists if}\ \int_{\R^2}\!\!g(x)dx=0,
\end{array}
\end{equation}
where $a(\e)$ is the one in (\ref{a}).
For the critical blow-up in low dimensions,
the situation is rather complicated because the rescaling argument is
no longer applicable.
Zhou \cite{Z92, Z93} proved that
there exist positive constants $c$ and $C$ independent of $\e$ 
(Hereafter in this section, we omit this description.) such that
\begin{equation}
\label{lifespan3}
\exp\left(c\e^{-p(p-1)}\right)\le T(\e)\le\exp\left(C\e^{-p(p-1)}\right)
\quad\mbox{for}\ p=p_0(n).
\end{equation}
\par
In higher dimensional case, $n\ge4$,
it is hard to get the same results as (\ref{lifespan1}) and (\ref{lifespan3})
because the fundamental solution is no longer positive.
Actually, we have
\begin{equation}
\label{lifespan4}
c\e^{-2p(p-1)/\gamma(p,n)+\sigma}\le T(\e)\le C\e^{-2p(p-1)/\gamma(p,n)}
\quad\mbox{for}\ 1<p<p_0(n),
\end{equation}
where $\sigma>0$ is a small error term.
The lower bound in (\ref{lifespan4}) was obtained by Di Pomponio and Georgiev \cite{DG01}.
On the other hand, the upper bound in (\ref{lifespan4})
is easily obtained by rescaling of the blowing-up solution in Sideris \cite{Si84}
which is stated in the history of Strauss' conjecture.
Such an argument can be found in Georgiev, Takamura and Zhou \cite{GTZ06}.
We note that it is possible to remove $\sigma$ in (\ref{lifespan4})
by assuming that the solution is radially symmetric.
See Section 6 in Lindblad and Sogge \cite{LS96}.
They also obtained the same lower bound as the one in (\ref{lifespan3}).
It is remarkable that, in $n=4$,
Li and Zhou \cite{LZ95} removed the assumption of radial symmetry
for the critical case as stated in the history on (\ref{GIVP}).
Their success depends on careful analysis in $L^2$ frame work.
Such a method is applicable to this case
because the nonlinear term is smooth by the fact that $p_0(4)=2$.
\par
As for the upper bound in (\ref{lifespan3}) for $n\ge4$,
following the proof in Zhou \cite{Z07} carefully,
one can find that $T(\e)\le \exp\left(\exp(C\e^{-p})\right)$.
Moreover, we point out that $T(\e)\le \exp(C\e^{-p^2})$
is implicitly obtained in Yordanov and Zhang \cite{YZ06}
if one follows their proof along with our argument.
See Remark \ref{rem:YZ} at the end of this paper.
But unfortunately both results are not optimal.
\par
In this paper, we prove the following expected theorem. 

\begin{thm}
\label{thm:main}
Let $n\ge4$ and $p=p_0(n)$.
Assume that both $f\in H^1(\R^n)$ and $g\in L^2(\R^n)$ are 
non-negative, do not vanish identically, and have compact support
such as $\{x\in\R^n\ :\ |x|\le R\}$, where $R$ is a positive constant.
Suppose that the problem {\rm(\ref{IVP})} has a solution
$(u,u_t)\in C([0,T(\e)),H^1(\R^n)\times L^2(\R^n))$ with
\begin{equation}
\label{support}
\mbox{\rm supp}(u,u_t)\subset\{(x,t)\in\R^n\times[0,\infty)\ :\ |x|\le t+R\}.
\end{equation}
Then, there exists a positive constant $\e_0
=\e_0(f,g,n,p,R)$ such that $T(\e)$ has to satisfy
\begin{equation}
\label{lifespan5}
T(\e)\le\exp\left(C\e^{-p(p-1)}\right)\quad\mbox{for}\ 0<\e\le\e_0,
\end{equation}
where $C$ is a positive constant independent of $\e$.
\end{thm}

\begin{rem}
\label{rem:regularity}
The differentiability of $\d\int_{\R^n}\!\!\!u(x,t)dx$ twice in $t$ follows from
the assumption on the regularity which is
the same as Yordanov and Zhang \cite{YZ06}.
See Sideris \cite{Si84} for details.
\end{rem}

Our success depends on the iteration argument of $L^p$ norm of the solution.
This is carried out on the integral inequality of the norm
which follows from $L^p$ boundedness of the maximal function via
Radon transform by Yordanov and Zhang \cite{YZ06}.
After repeating the estimates finitely many times
till $L^p$ norm is large enough,
we will be able to apply the blow-up theorem
for ordinary differential inequality
with the best condition only.

\section{Blow-up for ODI with a critical balance}
We shall start with the following blow-up result
for ordinary differential inequality.
This lemma is a modified version of Lemma 2.1 in Yordanov and Zhang \cite{YZ06}.
The key items are concrete expressions in (\ref{K_0T_1}) below.

\begin{lem}
\label{lem:ODI}
Let $p>1,\ a>0$ and $(p-1)a=q-2$. Suppose that $G\in C^2([0,T))$ satisfies
\begin{equation}
\label{ODI}
\left\{
\begin{array}{ll}
G(t)\ge Kt^a & \mbox{for}\ t\ge T_0,\\
G''(t)\ge B(t+R)^{-q}|G(t)|^p & \mbox{for}\ t\ge0,\\
G(0)>0,\quad G'(0)>0,
\end{array}
\right.
\end{equation}
where $B,K,R,T_0$ are positive constants with $T_0\ge R$.
Then, $T$ must satisfy that $T\le 2T_1$ provided $K\ge K_0$, where
\begin{equation}
\label{K_0T_1}
K_0=\left\{\frac{1}{2^{q/2}a}\sqrt{\frac{B}{p+1}}
\left(1-\frac{1}{2^{a\delta}}\right)\right\}^{-2/(p-1)},
\
T_1=\max\left\{T_0,\frac{G(0)}{G'(0)}\right\}
\end{equation}
with an arbitrarily chosen $\delta$ satisfying $0<\delta<(p-1)/2$.
\end{lem}

\par\noindent
{\bf Proof.} We prove this lemma by contradiction.
Assume that $T>2T_1$.
First we note that the second and third inequalities in (\ref{ODI}) yield
\begin{equation}
\label{growth}
G'(t)\ge G'(0)>0,\quad G(t)\ge G'(0)t+G(0)\ge G(0)>0\quad\mbox{for}\ t\ge0.
\end{equation}
Multiplying the second inequality in (\ref{ODI}) by $G'(t)$ and
integrating it over $[0,t]$, we have
\[
\begin{array}{ll}
\d\frac{1}{2}G'(t)^2
&\d\ge B\int_0^t(s+R)^{-q}G(s)^pG'(s)ds+\frac{1}{2}G'(0)^2\\
&\d>\frac{B}{(p+1)(t+R)^q}\left\{G(t)^{p+1}-G(0)^{p+1}\right\}\\
&\d\ge\frac{B}{(p+1)(t+R)^q}G(t)^p\left\{G(t)-G(0)\right\}
\end{array}
\]
for $t\ge0$. Restricting the time interval to $t\ge G(0)/G'(0)$
and making use of (\ref{growth}), we get
\[
\frac{1}{2}G(t)-G(0)\ge\frac{1}{2}\left\{G'(0)t-G(0)\right\}\ge0.
\]
Hence we obtain
\[
G'(t)>\sqrt{\frac{B}{p+1}}\cdot\frac{G(t)^{(p+1)/2}}{(t+R)^{q/2}}
\quad\mbox{for}\quad t\ge\frac{G(0)}{G'(0)}.
\]
\par
If $t\ge T_1(\ge R)$, one can make use of the first inequality in (\ref{ODI})
to obtain
\[
\frac{G'(t)}{G(t)^{1+\delta}}
>\sqrt{\frac{B}{p+1}}\cdot\frac{G(t)^{(p-1)/2-\delta}}{(t+R)^{q/2}}
\ge\sqrt{\frac{B}{p+1}}\cdot\frac{K^{(p-1)/2-\delta}}{2^{q/2}t^{q/2-a\{(p-1)/2-\delta\}}}
\]
for any $\delta$ satisfying $0<\delta<(p-1)/2$.
Noticing that $q/2-a(p-1)/2=1$
and integrating this inequality over $[T_1,t]$, we have
\[
\frac{1}{\delta}\left(\frac{1}{G(T_1)^\delta}-\frac{1}{G(t)^\delta}\right)
>\frac{1}{2^{q/2}a\delta}\sqrt{\frac{B}{p+1}}K^{(p-1)/2-\delta}
\left(\frac{1}{T_1^{a\delta}}-\frac{1}{t^{a\delta}}\right).
\]
Then, one can put $t=2T_1$ because of $T>2T_1$.
Neglecting $1/G(t)^\delta>0$ in the left hand side
and making use of the first inequality in (\ref{ODI}) with $t=T_1$, we obtain
\[
\frac{1}{K^\delta}\ge\left(\frac{T_1^a}{G(T_1)}\right)^\delta
>\frac{1}{2^{q/2}a}\sqrt{\frac{B}{p+1}}\left(1-\frac{1}{2^{a\delta}}\right)
K^{(p-1)/2-\delta}.
\]
This inequality contradicts to the choice of $K\ge K_0$.
Therefore we conclude that $T\le 2T_1$. The lemma is now established.
\hfill$\Box$

\section{Growing up of $L^p$ norm of the solution}
In this section, we shall construct an iteration of estimates
for $L^p$ norm of the solution.
As stated in Remark \ref{rem:regularity},
the assumption on the regularity in Theorem \ref{thm:main} yields
\[
F(t)=\int_{\R^n}\!\!\!u(x,t)dx\in C^2([0,T(\e)),
\]
so that we have
\[
F''(t)=\int_{\R^n}|u(x,t)|^pdx=\|u(\cdot,t)\|_{L^p(\R^n)}^p.
\]
\par
The iteration argument will give us an enough growth of the norm for large time.
To this end, we have to start with the following basic frame of the iteration.

\begin{prop}
\label{prop:frame}
Suppose that the assumption in Theorem \ref{thm:main} is fulfilled.
Then, there exists a positive constant $C=C(f,g,n,p,R)$
such that $F(t)=\d\int_{\R^n}\!\!\!u(x,t)dx$ for $t\ge R$ satisfies
\begin{equation}
\label{frame}
F''(t)\ge C\int_0^{t-R}\frac{\rho^{(n-1)(1-p/2)}d\rho}{(t-\rho+R)^{(n-1)p/2}}
\left(\int_0^{(t-\rho-R)/2}F''(s)ds\right)^p.
\end{equation}
\end{prop}

\par\noindent
{\bf Proof.} This proposition immediately follows from the combination of
two estimates for Radon transformation,
(2.14) and (2.21), in Yordanov and Zhang \cite{YZ06}.
\hfill$\Box$
\vskip10pt
\par
The next proposition is the basic estimate for the first step of our iteration.

\begin{prop}
\label{prop:step0}
Suppose that the assumption in Theorem \ref{thm:main} is fulfilled.
Then, there exists a positive constant $C=C(f,g,n,p,R)$
such that $F(t)=\d\int_{\R^n}\!\!\!u(x,t)dx$ for $t\ge0$ satisfies
\begin{equation}
\label{step0}
F''(t)\ge C\e^p(t+R)^{(n-1)(1-p/2)}.
\end{equation}
\end{prop}

\par\noindent
{\bf Proof.} This is exactly (2.5') in Yordanov and Zhang \cite{YZ06}.
They employed a special test function.
Without such a technique,
the easy proof for slightly different data can be found in Rammaha \cite{R88},
in which the short and simple proof of Sideris' blow-up theorem in high dimensions
is given.
\hfill$\Box$

\begin{rem}
It is trivial that we can write the same $C$ in Propositions \ref{prop:frame} and
\ref{prop:step0}.
\end{rem}

\par
The main estimate in our iteration is the  following proposition.

\begin{prop}
\label{prop:stepj}
Suppose that the assumption in Theorem \ref{thm:main} is fulfilled.
Then, $F(t)=\d\int_{\R^n}\!\!\!u(x,t)dx$ for $t\ge a_jR\ (j=1,2,3,\cdots)$ satisfies that
\begin{equation}
\label{stepj}
F''(t)\ge C_j(t-a_jR)^{(n-1)(1-p/2)}
\left(\log\frac{t+(a_j-2)R}{2(a_j-1)R}\right)^{(p^j-1)/(p-1)}.
\end{equation}
Here we set $a_j=3\cdot4^{j-1}-1$ and
\begin{equation}
\label{C_j}
\begin{array}{l}
\d C_j=\exp\left\{p^{j-1}\left(\log(C_0C_1C_p^{-S(j)})\right)-\log C_0\right\}
\quad (j\ge2),\\
\d C_1=\frac{C^{p+1}}{2^{n-2}\cdot3^{(n-1)p/2}\{n-(n-1)p/2\}^p}\e^{p^2},
\end{array}
\end{equation}
where $C$ is the one in Propositions \ref{prop:frame}, \ref{prop:step0},
and
\begin{equation}
\label{notation}
C_0=\left\{\frac{(p-1)C}{2^{n-1+(n+1)p/2}\cdot3^{np-1}p}\right\}^{1/(p-1)},
\ C_p=2^{(n+1)p}p,
\ S(j)=\sum_{k=1}^{j-1}\frac{k}{p^k}. 
\end{equation}
\end{prop}

\par\noindent
{\bf Proof.} Recall that $1<p=p_0(n)\le p_0(4)=2$ for $n\ge4$.
First we shall show this proposition for $j=1$.
Replacing $F''(s)$ in the right hand side of (\ref{frame})
by the lower bound of $F''(t)$ in (\ref{step0}), we have
\[
F''(t)\ge C^{p+1}\e^{p^2}\int_0^{t-R}\frac{\rho^{(n-1)(1-p/2)}d\rho}{(t-\rho+R)^{(n-1)p/2}}
\left(\int_0^{(t-\rho-R)/2}s^{(n-1)(1-p/2)}ds\right)^p
\]
for $t\ge R$.
Hence it follows that
\[
F''(t)\ge \frac{2^{n-2}\cdot3^{(n-1)p/2}C_1}{2^{np-(n-1)p^2/2}}\int_0^{t-R}
\frac{\rho^{(n-1)(1-p/2)}(t-\rho-R)^{np-(n-1)p^2/2}}{(t-\rho+R)^{(n-1)p/2}}d\rho
\]
for $t\ge R$, where $C_1$ is defined in (\ref{C_j}).
From now on, we restrict the time interval to $t\ge a_1R=2R$
and diminish the domain of the $\rho$-integral to $[0,t-2R]$.
Then we have $t-\rho\ge 2R$ in the $\rho$-integral.
We now employ the following elementary lemma.

\begin{lem}
\label{lem:ineq}
Let $M$ and $R$ be positive constants.
Then $t-\rho\ge MR$ is equivalent to
\[
(M+1)\{t-\rho-(M-1)R\}\ge t-\rho+R.
\]
\end{lem}
It is easy to prove this lemma.
We omit the proof.
\par
Making use of Lemma \ref{lem:ineq} with $M=2$ and
the relation
\begin{equation}
\label{relation}
np-\frac{n-1}{2}p^2=\frac{n-1}{2}p-1
\end{equation}
which is equivalent to (\ref{gamma}), we obtain
\[
F''(t)\ge \frac{2^{n-2}C_1}{2^{(n-1)p/2-1}}\int_0^{t-2R}
\frac{\rho^{(n-1)(1-p/2)}}{t-\rho-R}d\rho
\]
for $t\ge 2R$.
Hence, cutting the domain of the $\rho$-integral to be an upper half,
we have
\[
\begin{array}{ll}
F''(t)
&\d \ge C_1(t-2R)^{(n-1)(1-p/2)}\int_{(t-2R)/2}^{t-2R}
\frac{1}{t-\rho-R}d\rho\\
&\d \ge C_1(t-2R)^{(n-1)(1-p/2)}\log\frac{t}{2R}
\end{array}
\]
for $t\ge 2R$.
Therefore (\ref{stepj}) is true for $j=1$.
\par
Next we shall show (\ref{stepj}) by induction.
Assume that (\ref{stepj}) holds but $C_j$ is unknown except for $j=1$.
Later we look for the relation between $C_j$ and $C_{j+1}$
which yields (\ref{C_j}).
To this end, we restrict the time interval $t\ge a_jR$ to $t\ge(2a_j+1)R$.
Then it follows from (\ref{frame}) that
\[
F''(t)\ge C\int_0^{t-(2a_j+1)R}\frac{\rho^{(n-1)(1-p/2)}d\rho}{(t-\rho+R)^{(n-1)p/2}}
\left(\int_{a_jR}^{(t-\rho-R)/2}F''(s)ds\right)^p
\]
for $t\ge(2a_j+1)R$.
Making use of (\ref{stepj}), we have
\[
F''(t)\ge CC_j^p\int_0^{t-(2a_j+1)R}
\frac{\rho^{(n-1)(1-p/2)}\{I_j(t,\rho)\}^p}{(t-\rho+R)^{(n-1)p/2}}d\rho
\]
for $t\ge(2a_j+1)R$, where we set
\[
I_j(t,\rho)=\int_{a_jR}^{(t-\rho-R)/2}\hspace{-30pt}(s-a_jR)^{(n-1)(1-p/2)}
\left(\log\frac{s+(a_j-2)R}{2(a_j-1)R}\right)^{(p^j-1)/(p-1)}ds.
\]
\par
Now we restrict the time interval further to $t\ge 2(a_j+1)R$
and diminish the domain of the $\rho$-integral to $[0,t-2(a_j+1)R]$.
Then we have $t-\rho\ge 2(a_j+1)R$ in the $\rho$-integral.
We note that one can diminish also the domain of the $s$-integral
to $[a_j(t-\rho-R)/(2a_j+1),(t-\rho-R)/2]$ because of
\[
a_jR\le\frac{a_j}{2a_j+1}(t-\rho-R).
\]
Since $(s-a_jR)$ in the $s$-integral is estimated by
\[
\frac{a_j}{2a_j+1}(t-\rho-R)-a_jR=\frac{a_j}{2a_j+1}(t-\rho-2(a_j+1)R),
\]
we have
\[
\begin{array}{ll}
I_j(t,\rho)\ge
&\d \left(\frac{t-\rho-2(a_j+1)R}{3}\right)^{(n-1)(1-p/2)}\times\\
&\d \times\int_{(t-\rho-R)a_j/(2a_j+1)}^{(t-\rho-R)/2}
\left(\log\frac{s+(a_j-2)R}{2(a_j-1)R}\right)^{(p^j-1)/(p-1)}ds.
\end{array}
\]
Moreover, it follows from $t-\rho\ge 2(a_j+1)R$ that
the variable in the logarithmic term is estimated as
\[
\begin{array}{l}
\d \frac{(t-\rho-R)a_j/(2a_j+1)+(a_j-2)R}{2(a_j-1)R}\\
\d =\frac{t-\rho+(2a_j-1)R}{2(2a_j+1)R}
+\frac{t-\rho-(a_j+3)R}{2(a_j-1)(2a_j+1)R}\\
\d \ge\frac{t-\rho+(2a_j-1)R}{2(2a_j+1)R}
+\frac{(a_j-1)R}{2(a_j-1)(2a_j+1)R}.
\end{array}
\]
Hence, neglecting the last positive term in the above inequality, we get
\[
\begin{array}{ll}
I_j(t,\rho)\ge
&\d \left(\frac{t-\rho-2(a_j+1)R}{3}\right)^{(n-1)(1-p/2)}\times\\
&\d \times\frac{t-\rho-R}{2(2a_j+1)}
\left(\log\frac{t-\rho+(2a_j-1)R}{2(2a_j+1)R}\right)^{(p^j-1)/(p-1)}.
\end{array}
\]
Therefore (\ref{relation}) yields
\[
\begin{array}{l}
\d F''(t)\ge\frac{CC_j^p}{3^{(n-1)p/2-1}(2a_j)^p}\int_0^{t-2(a_j+1)R}
\frac{\rho^{(n-1)(1-p/2)}}{(t-\rho+R)^{(n-1)p/2}}d\rho\times\\
\d \times\{t-\rho-2(a_j+1)R\}^{(n-1)p/2-1}
\left(\log\frac{t-\rho+(2a_j-1)R}{2(2a_j+1)R}\right)^{(p^{j+1}-p)/(p-1)}
\end{array}
\]
for $t\ge 2(a_j+1)R$.
\par
Now we restrict the time interval again to $t\ge(2a_j+3)R$.
Then it follows from Lemma \ref{lem:ineq} with $M=2a_j+3$ that
\[
t-\rho+R\le(2a_j+4)\{t-\rho-2(a_j+1)R\}\le2^2a_j\{t-\rho-2(a_j+1)R\}.
\]
Hence we have
\[
\begin{array}{l}
\d F''(t)\ge\frac{CC_j^p}{2^{np}\cdot3^{(n-1)p/2-1}a_j^{(n+1)p/2}}
\times\\
\d \times
\int_0^{t-(2a_j+3)R}\hspace{-20pt}
\frac{\rho^{(n-1)(1-p/2)}}{t-\rho-2(a_j+1)R}
\left(\log\frac{t-\rho+(2a_j-1)R}{2(2a_j+1)R}\right)^{(p^{j+1}-p)/(p-1)}d\rho
\end{array}
\]
for $t\ge(2a_j+3)R$.
This inequality implies
\[
\begin{array}{ll}
F''(t)\ge
&\d \frac{CC_j^p}{2^{n-1+(n+1)p/2}\cdot3^{(n-1)p/2-1}a_j^{(n+1)p/2}}
\{t-(2a_j+3)R\}^{(n-1)(1-p/2)}\times\\
&\d \times
\int_{\{t-(2a_j+3)R\}/2}^{t-(2a_j+3)R}
\frac{\left(\log\frac{t-\rho+(2a_j-1)R}{2(2a_j+1)R}\right)^{(p^{j+1}-p)/(p-1)}}{t-\rho+(2a_j-1)R}
d\rho
\end{array}
\]
for $t\ge(2a_j+3)R$.
Noticing that
\[
\frac{p^{j+1}-p}{p-1}+1=\frac{p^{j+1}-1}{p-1}\le\frac{p^{j+1}}{p-1}
\]
and
\[
a_j=3\cdot4^{j-1}-1\le3\cdot2^{2j},
\]
we obtain
\[
\begin{array}{ll}
\d F''(t)\ge
&\d \frac{CC_j^p}{2^{n-1+(n+1)p/2}\cdot3^{np-1}\cdot2^{(n+1)pj}}
\{t-(2a_j+3)R\}^{(n-1)(1-p/2)}\times\\
&\d \times
\frac{p-1}{p^{j+1}}
\left(\log\frac{t+(6a_j+1)R}{2(4a_j+2)R}\right)^{(p^{j+1}-1)/(p-1)}
\end{array}
\]
for $t\ge(2a_j+3)R$.
Therefore it follows from $a_{j+1}=4a_j+3$ that
\[
\begin{array}{ll}
\d F''(t)\ge
&\d \frac{(p-1)CC_j^p}{2^{n-1+(n+1)p/2}\cdot3^{np-1}p\cdot(2^{(n+1)p}p)^j}
\times\\
&\d \times
(t-a_{j+1}R)^{(n-1)(1-p/2)}
\left(\log\frac{t+(a_{j+1}-2)R}{2(a_{j+1}-1)R}\right)^{(p^{j+1}-1)/(p-1)}
\end{array}
\]
for $t\ge a_{j+1}R$.
\par
As a conclusion, if $C_j$ is defined by
\[
C_{j+1}=\frac{C_0^{p-1}C_j^p}{C_p^j}
\quad(j\ge1),
\]
where $C_0$ and $C_p$ are defined by (\ref{notation}),
then (\ref{stepj}) is valid for all $j\ge1$.
This equality is rewritten as
\[
\log C_{j+1}=p\log C_j-j\log C_p+\log C_0^{p-1}.
\]
It is clear that $C_2$ defined by this equality is the one in (\ref{C_j}).
For $j\ge2$, we have the following concrete expression of $\log C_{j+1}$ inductively. 
\[
\begin{array}{ll}
\log C_{j+1}
&\d =p^j\log C_1-\sum_{k=1}^jkp^{j-k}\log C_p+\sum_{k=0}^{j-1}p^k\log C_0^{p-1}\\
&\d =p^j\left\{\log C_1-S(j+1)\log C_p+\log C_0\right\}-\log C_0.
\end{array}
\]
This is exactly (\ref{C_j}).
Therefore Proposition \ref{prop:stepj} is now established.
\hfill$\Box$


\section{Upper bound of the lifespan}
In this section, we complete the proof of Theorem \ref{thm:main}.
The first step is to shift the estimate for $F''(t)=\|u(\cdot,t)\|^p_{L^p(\R^n)}$
to the one for $F(t)=\d\int_{\R^n}\!\!\!u(x,t)dx$.
One of the key in this section is the assumption on the initial data in Theorem \ref{thm:main},
\begin{equation}
\label{F(0)}
F(0)=\e\int_{\R^n}\!\!f(x)dx>0,\quad F'(0)=\e\int_{\R^n}\!\!g(x)dx>0.
\end{equation}
This yields that $F(t)>0$ and $F'(t)>0$ for $t\ge0$.
Because it follows from the support condition (\ref{support}) and H\"older's inequality
that
\begin{equation}
\label{F''}
F''(t)\ge\left\{\mbox{vol}({\bf B}^n(0,1))\right\}^{1-p}(t+R)^{-n(p-1)}|F(t)|^p
\end{equation}
for $t\ge0$, where $\mbox{vol}({\bf B}^n(0,1))$ is a volume of a unit ball in $\R^n$.
\par
Now we start with the following proposition.

\begin{prop}
\label{prop:F}
Suppose that the assumption in Theorem \ref{thm:main} is fulfilled.
Then, $F(t)=\d\int_{\R^n}\!\!\!u(x,t)dx$ for $t\ge \{2(a_j+2)R\}^2\ (j=1,2,3,\cdots)$
satisfies
\begin{equation}
\label{F}
F(t)\ge\frac{C_j}{16^jD}
\left(\frac{1}{2}\log t\right)^{(p^j-1)/(p-1)}t^{n+1-(n-1)p/2},
\end{equation}
where $D=3^2\cdot2^{3n-2-3(n-1)p/2}$,
$a_j$ and $C_j$ are defined in Proposition \ref{prop:stepj}.
\end{prop}

\par\noindent
{\bf Proof.} 
Integrating (\ref{stepj}) in Proposition \ref{prop:stepj}
over $[a_jR,t]$, we have
\[
F'(t)\ge C_j\int_{a_jR}^t(s-a_jR)^{(n-1)(1-p/2)}
\left(\log\frac{s+(a_j-2)R}{2(a_j-1)R}\right)^{(p^j-1)/(p-1)}ds
\]
for $t\ge a_jR$.
Here we restrict the time interval to $t\ge(a_j+1)R$
and diminish the domain of the $s$-integral to $[a_jt/(a_j+1),t]$.
$(s-a_jR)$ in the integral is estimated by
\[
\frac{a_j}{a_j+1}t-a_jR\ge\frac{1}{2}\{t-(a_j+1)R\}.
\]
Also the variable of the logarithmic term is estimated by
\[
\frac{a_jt/(a_j+1)+(a_j-2)R}{2(a_j-1)R}
=\frac{t+(a_j+1)R}{2(a_j+1)R}
+\frac{t-(a_j+1)R}{2(a_j-1)(a_j+1)R}.
\]
Hence we obtain
\[
F'(t)\ge\frac{C_j\{t-(a_j+1)R\}^{n-(n-1)p/2}}{2^{n-(n-1)p/2}a_j}
\left(\log\frac{t+(a_j+1)R}{2(a_j+1)R}\right)^{(p^j-1)/(p-1)}
\]
for $t\ge (a_j+1)R$.
\par
Integrating this inequality over $[(a_j+1)R,t]$, we have
\[
\begin{array}{ll}
\d F(t)\ge\frac{C_j}{2^{n-(n-1)p/2}a_j}
&\d \int_{(a_j+1)R}^t\{s-(a_j+1)R\}^{n-(n-1)p/2}\times\\
&\d \times\left(\log\frac{s+(a_j+1)R}{2(a_j+1)R}\right)^{(p^j-1)/(p-1)}ds
\end{array}
\]
for $t\ge(a_j+1)R$.
Similary to the above, we restrict the time interval to $t\ge(a_j+2)R$
and diminish the domain of the $s$-integral to $[(a_j+1)t/(a_j+2),t]$.
$(s-(a_j+1)R)$ in the integral is estimated by
\[
\frac{a_j+1}{a_j+2}t-(a_j+1)R\ge\frac{1}{2}\{t-(a_j+2)R\}.
\]
Also the variable of the logarithmic term is estimated by
\[
\frac{(a_j+1)t/(a_j+2)+(a_j+1)R}{2(a_j+1)R}
=\frac{t+(a_j+2)R}{2(a_j+2)R}.
\]
Hence we obtain
\[
F(t)\ge\frac{C_j\{t-(a_j+2)R\}^{n+1-(n-1)p/2}}{2^{2n+1-(n-1)p}a_j^2}
\left(\log\frac{t+(a_j+2)R}{2(a_j+2)R}\right)^{(p^j-1)/(p-1)}
\]
for $t\ge(a_j+2)R$.
\par
Restricting the time interval further to $t\ge2(a_j+2)R$, we have
\[
F(t)\ge\frac{C_jt^{n+1-(n-1)p/2}}{2^{3n+2-3(n-1)p/2}a_j^2}
\left(\log\frac{t}{2(a_j+2)R}\right)^{(p^j-1)/(p-1)}.
\]
Note that we may assume $2(a_1+2)R\ge1$ without loss of the generality.
Therefore we finally obtain
\[
F(t)\ge\frac{C_jt^{n+1-(n-1)p/2}}{2^{3n+2-3(n-1)p/2}a_j^2}
\left(\frac{1}{2}\log t\right)^{(p^j-1)/(p-1)}
\]
for $t\ge\{2(a_j+2)R\}^2\ge2(a_j+2)R$.
The proof is now ended by trivial inequality $a_j\le 3\cdot2^{-2}\cdot4^j$.
\hfill$\Box$

\vskip10pt
\par\noindent
{\bf Proof of Theorem \ref{thm:main}.}
Let $j\ge2$.
Define a sequence of time interval $\{I(j)\}$ by
\begin{equation}
\label{I}
I(j)=\left[\{2(a_j+2)R\}^2,\{2(a_{j+1}+2)R\}^2\right]
\end{equation}
and set
\[
K_j(t)=\frac{C_j}{16^jD}\left(\frac{1}{2}\log t\right)^{(p^j-1)/(p-1)}
\]
which is the coefficient of $t^{n+1-(n-1)p/2}$ in (\ref{F}).
Then it follows from the definition of $C_j$ in (\ref{C_j}) that
\[
K_j(t)=\exp\left\{p^{j-1}\log L_j(t)
-j\log16-\log(C_0D)-\frac{\log\left(\log \sqrt{t}\right)}{p-1}\right\},
\]
where we set
\[
L_j(t)=C_0C_1C_p^{-S(j)}\left(\frac{1}{2}\log t\right)^{p/(p-1)}.
\]
In view of the definition of $C_1$ in (\ref{C_j}),
we have $L_j(t)\ge e$
provided
\begin{equation}
\label{lifespan}
\e^{p(p-1)}\log t\ge E,
\end{equation}
where
\[
E=2\left(\frac{2^{n-2}\cdot3^{(n-1)p/2}\{n-(n-1)p/2\}^p\cdot eC_p^{S(\infty)}}
{C_0C^{p+1}}\right)^{(p-1)/p}>0.
\]
Because $S(j)$ is monotonously increasing in $j$,
but converges to a positive constant $S(\infty)$. 
\par
From now on, we assume (\ref{lifespan}).
Then it follows that
\[
K_j(t)\ge\exp\left\{p^{j-1}
-j\log16-\log(C_0D)-\frac{\log\left(\log\{2(a_{j+1}+2)R\}\right)}{p-1}\right\}
\]
for $t\in I(j)$.
We note that the right hand side of this inequality
goes to infinity if $j$ tends to infinity.
Hence, for $K_0$ defined in (\ref{K_0T_1}) with $a=n+1-(n-1)p/2>0$ and
$B=\left\{\mbox{vol}({\bf B}^n(0,1))\right\}^{1-p}>0$,
there exists an integer $J=J(f,g,n,p,R)$ such that
\[
F(t)\ge K_0t^{n+1-(n-1)p/2}\quad\mbox{for}\ t\in I(j)
\]
as far as $j\ge J$. Therefore the definition of $I(j)$ implies
\[
F(t)\ge K_0t^{n+1-(n-1)p/2}\quad\mbox{for}\ t\ge\{2(a_J+2)R\}^2.
\]
\par
Now we are in a position to apply Lemma \ref{lem:ODI} to our situation with
\[
G=F,\ B=\left\{\mbox{vol}({\bf B}^n(0,1))\right\}^{1-p}
\]
and
\[
a=n+1-\frac{n-1}{2}p,\ q=n(p-1)
\]
because of (\ref{F''}).
We note that the condition $(p-1)a=q-2$ in this setting is equivalent to $p=p_0(n)$.
First we set
\[
T_0(\e)=\exp\left(E\e^{-p(p-1)}\right),
\]
where $E$ is the one in (\ref{lifespan}).
Then there exists $\e_0=\e_0(f,g,n,p,R)$ such that
\[
T_0(\e)\ge\{2(a_J+2)R\}^2\quad\mbox{and}\quad
2\max\left\{T_0(\e),\frac{F(0)}{F'(0)}\right\}\le\exp\left(2E\e^{-p(p-1)}\right)
\]
hold for $0<\e\le\e_0$
because $J$ and $F(0)/F'(0)$ are independent of $\e$ as we see.
If the lifespan $T(\e)$ satisfies $T(\e)>T_0(\e)$, then we have
\[
F(t)\ge K_0 t^{n+1-(n-1)p/2}\quad\mbox{for}\ t\in[T_0(\e),T(\e))
\]
by definition of $T_0(\e)$
because such a $t$ satisfies $\e^{p(p-1)}\log t\ge E$.
Lemma \ref{lem:ODI} says that this inequality implies
\[
t\le2\max\left\{T_0(\e),\frac{F(0)}{F'(0)}\right\}\le\exp\left(2E\e^{-p(p-1)}\right).
\]
Taking a supremum over $t\in[T_0(\e),T(\e))$,
we get
\begin{equation}
\label{lifespanfinal}
T(\e)\le\exp\left(2E\e^{-p(p-1)}\right)
\quad\mbox{for}\ 0<\e\le\e_0.
\end{equation}
The counter case $T(\e)\le T_0(\e)$ is trivial.
Therefore (\ref{lifespanfinal}) holds for any cases.
The proof of Theorem \ref{thm:main} is now completed.
\hfill$\Box$

\begin{rem}
\label{rem:YZ}
It is easy to check that
the blow-up condition in Yordanov and Zhang \cite{YZ06} is
\[
\lim_{t\rightarrow\infty}\e^{p^2}\log t=\infty.
\]
But one can find that their estimate is equivalent to
Proposition \ref{prop:stepj} with $j=1$.
Hence, applying the above argument to such an estimate,
we have
\[
T(\e)\le\exp\left(2\wt{E}\e^{-p^2}\right)\quad\mbox{for}\quad0<\e\le\e_0
\]
with a different constant $\wt{E}>0$ from $E$.
This result is stated in Introduction.
The improvement of the upper bound of the lifespan is carried out by our iteration argument. 
\end{rem}

\bibliographystyle{plain}

\end{document}